\documentclass[12pt]{amsart}
\usepackage{graphicx, amssymb, colordvi,amsmath}
\usepackage{psfrag}
\usepackage{pstricks}
\usepackage[bookmarks=false,colorlinks,linkcolor=mylinkcolor,citecolor=mycitecolor]{hyperref}
   \definecolor{mycitecolor}{rgb}{.2,.6,.2}
   \definecolor{mylinkcolor}{rgb}{.2,.2,.6}

\newtheorem{theorem}{Theorem}
\newtheorem{conjecture}{Conjecture}

\newtheorem{lemma}[theorem]{Lemma}
\newtheorem{remark}[theorem]{Remark}
\newcommand{\harpoon}{\rightharpoonup}
\newcommand{\Sh}{\mbox{\rm sh}}
\newcommand{\sh}{\mbox{\scriptsize\rm sh}}
\newcommand{\demph}[1]{\Blue{\emph{#1}}}

\newcommand{\eqsj}{\equiv_{\textsc{sj}}}
\newcommand{\eqsp}{\equiv_{\textsc{sp}}}
\newcommand{\ljoin}{{\triangleleft}}
\newcommand{\ujoin}{\hbox{\small${\vartriangle}$}}
\newcommand{\join}{\lozenge}
\renewcommand{\ljoin}{\mathchoice%
  {\raisebox{.17\height}{\small$\displaystyle\triangleleft$}}
  {\raisebox{.17\height}{\small$\textstyle\triangleleft$}}
  {\raisebox{.17\height}{\small$\scriptstyle\triangleleft$}}
  {\raisebox{.17\height}{\small$\scriptscriptstyle\triangleleft$}}
}
\renewcommand{\ujoin}{\mathchoice%
  {\raisebox{.19\height}{\scriptsize$\displaystyle\vartriangle$}}
  {\raisebox{.19\height}{\scriptsize$\textstyle\vartriangle$}}
  {\raisebox{.19\height}{\scriptsize$\scriptstyle\vartriangle$}}
  {\raisebox{.19\height}{\scriptsize$\scriptscriptstyle\vartriangle$}}
}
\renewcommand{\join}{\mathchoice%
  {\raisebox{.14\height}{\footnotesize$\displaystyle\lozenge$}}
  {\raisebox{.14\height}{\footnotesize$\textstyle\lozenge$}}
  {\raisebox{.14\height}{\footnotesize$\scriptstyle\lozenge$}}
  {\raisebox{.14\height}{\footnotesize$\scriptscriptstyle\lozenge$}}
}
\newcommand{\des}{\mathsf d}
\newcommand{\word}{\mathsf w}
\newcommand{\tild}[1]{#1\tilde{\phantom{e}}}
\renewcommand{\tild}[1]{{#1\,}'}
\newcommand{\warrow}[1]{{\,\rightsquigarrow_{#1}\,}}
\newcommand{\sarrow}[1]{{\,\rightarrow_{#1}\,}}
\newcommand{\C}{{\mathbb C}}
\newcommand{\Gr}{{\mathit{Gr}}}
\newcommand{\Z}{\mathcal{Z}}

\copyrightinfo{}{}
\date{29 October 2009}
\title[Skew Littlewood-Richardson rules from Hopf algebras]%
  {Skew Littlewood-Richardson rules\\ from Hopf algebras}

\author{Thomas Lam}
\address[Lam]{
    Department of Mathematics\\
         University of Michigan\\
         Ann Arbor, MI \ 48109
         }
\email{tfylam@umich.edu}
  \urladdr{http://www.math.lsa.umich.edu/\~{}tfylam/}

\author{Aaron Lauve}
\address[Lauve]{
    Department of Mathematics\\
         Texas A\&M University\\
         College Station, TX \ 77843
         }
\email{lauve@math.tamu.edu}
  \urladdr{http://www.math.tamu.edu/\~{}lauve}

\author{Frank Sottile}
\address[Sottile]{Department of Mathematics\\
         Texas A\&M University\\
         College Station\\
         TX \ 77843}
\email{sottile@math.tamu.edu}
\urladdr{http://www.math.tamu.edu/\~{}sottile/}
\thanks{Sottile was supported by the NSF grant DMS-0701050.}
\thanks{Lam was supported by a Sloan Fellowship, and by NSF grants DMS-0901111 and
DMS-0652641.}


\begin{document}

\begin{abstract}
 We use Hopf algebras to prove a version
 of the Littlewood-Richardson rule for skew Schur functions,
 which implies a conjecture of Assaf and McNamara.
 We also establish skew Littlewood-Richardson rules for Schur $P$- and
 $Q$-functions and noncommutative ribbon Schur functions, as well as skew Pieri rules for
 $k$-Schur functions, dual $k$-Schur functions, and for the homology of the affine Grassmannian
 of the symplectic group.
\end{abstract}
\maketitle

Assaf and McNamara~\cite{AM} recently used combinatorics to give an
elegant and surprising formula for the product of a skew Schur
function and a complete homogeneous symmetric function.  Their paper
included a conjectural skew version of the Littlewood-Richardson
rule, and also an appendix by one of us (Lam) with a simple
algebraic proof of their formula.  We show how these formulas and
much more are special cases of a simple formula that holds for any
pair of dual Hopf algebras. We first establish this Hopf-algebraic
formula, and then apply it to obtain formulas in some well-known
Hopf algebras in combinatorics.

%
\section{A Hopf algebraic formula}\label{sec: harpoon}

We assume basic familiarity with Hopf algebras, as found in the opening chapters of
the book~\cite{Mont}.
Let \Blue{$H$}, \Blue{$H^*$} be a pair of dual Hopf algebras over a field $\Bbbk$.
This means that there is a  nondegenerate pairing
$\langle\cdot,\cdot\rangle\colon H\otimes H^*\to \Bbbk$ for which the structure of $H^*$ is
dual to that of $H$ and vice-versa.
For example, $H$ could be finite-dimensional and $H^*$ its linear dual, or $H$ could be graded
with each component finite-dimensional and $H^*$ its graded dual.
These algebras naturally act on each other~\cite[1.6.5]{Mont}:
suppose that $h\in H$ and $a\in H^*$ and set
 \begin{equation}\label{Eq:harpoon}
   \Blue{h\harpoon a}\ :=\ \sum \langle h, a_2\rangle a_1
   \qquad\mbox{and}\qquad
   \Blue{a\harpoon h}\ :=\ \sum \langle h_2, a\rangle h_1\,.
 \end{equation}
(We use Sweedler notation for the coproduct,
 $\Delta h \ =\ \sum h_1\otimes h_2$.)
These left actions are the adjoints of right multiplication:
for $g,h\in H$ and $a,b\in H^*$,
\[
   \langle g, h\harpoon a\rangle\ =\ \langle g\cdot h, a \rangle
   \qquad\mbox{and}\qquad
   \langle a\harpoon h, b\rangle\ =\ \langle h, b\cdot a \rangle\,.
\]
This shows that $H^*$ is a left $H$-module under the action in \eqref{Eq:harpoon}.
In fact, $H^*$ is a left $H$--module algebra, meaning that for $a,b\in H^*$ and $h\in H$,
 \begin{equation}\label{Eq:CoModAlg}
   h\harpoon(a\cdot b)\ =\ \sum (h_1\harpoon a)\cdot(h_2\harpoon b)\,.
 \end{equation}

Recall that the \demph{counit} $\varepsilon\colon H \to \Bbbk$ and \demph{antipode}
$S\colon H \to H$ satisfy $\sum h_1 \cdot \varepsilon(h_2) = h$ and
$\sum h_1\cdot S(h_2)=\varepsilon(h)\cdot 1_H$ for all $h\in H$.

\begin{lemma}\label{Lem:Hopf}
 For $g,h\in H$ and $a\in H^*$, we have
 \begin{equation}\label{Two}
   (a\harpoon g)\cdot h \ =\  \sum (S(h_2)\harpoon a)\harpoon (g\cdot h_1)\,.
 \end{equation}
\end{lemma}

\begin{proof}
 Let $b\in H^*$.
 We prove first the formula
 \begin{equation}\label{One}
    (h\harpoon b)\cdot a \ =\ \sum h_1\harpoon (b\cdot (S(h_2)\harpoon a))\,.
 \end{equation}
 (This is essentially $(*)$ in the proof of Lemma 2.1.4 in~\cite{Mont}.)
 Expanding the sum using~\eqref{Eq:CoModAlg} and coassociativity,
 $(\Delta\otimes 1) \circ \Delta(h)=(1\otimes\Delta)\circ \Delta(h)=
   \sum h_1\otimes h_2\otimes h_3$, gives
 \begin{align}
   \sum h_1\harpoon (b\cdot (S(h_2)\harpoon a)) \ &=\
    \sum (h_1\harpoon b)\cdot (h_2\harpoon (S(h_3)\harpoon a))\nonumber\\
   \ &=\ \sum (h_1\harpoon b)\cdot ((h_2\cdot S(h_3))\harpoon a)\label{first}\\
  \  &=\ (h\harpoon b)\cdot a \,.   \label{second}
 \end{align}
Here, \eqref{first} follows as $H^*$ is an $H$-module and~\eqref{second} from
 the antipode and counit conditions.

 Note that
  $\langle (a\harpoon g)\cdot h, b \rangle =
  \langle a\harpoon g, h \harpoon b \rangle =
  \langle g, (h\harpoon b)\cdot a\rangle$.
 Using~\eqref{One} this becomes
 \begin{align*}
   \bigl\langle g, \sum h_1\harpoon (b\cdot (S(h_2)\harpoon a))\bigr\rangle
   \ &=\ \sum \bigl\langle g\cdot h_1, b\cdot ( S(h_2)\harpoon a)\bigr\rangle\\
   \ &=\ \bigl\langle \;\sum (S(h_2)\harpoon a) \harpoon (g\cdot h_1), b\bigr\rangle\, ,
 \end{align*}
 which proves the lemma, as this holds for all $b\in H^*$.
\end{proof}

\begin{remark}
{\rm
 This proof is identical to the argument in the appendix to~\cite{AM},
where $h$ was a complete homogeneous symmetric function in the
 Hopf algebra $H$ of symmetric functions.
}
\end{remark}

\section{Application to distinguished bases}\label{sec: bases}

We apply Lemma \ref{Lem:Hopf} to produce skew Littlewood-Richardson rules for several Hopf algebras in algebraic combinatorics. We isolate the common features of those arguments.

In the notation of Section \ref{sec: harpoon}, let $\{L_\lambda\}\subset H$ and $\{R_\lambda\}\subset H^*$ be dual bases indexed by some set $\mathcal P$, so $\langle L_\lambda, R_\mu \rangle = \delta_{\lambda, \mu}$ for $\lambda,\mu\in\mathcal P$. Define structure constants for $H$ and $H^*$ via
\begin{align}
\label{eq: L structure}
    L_\lambda \cdot L_\mu \ &=\ \sum_{\nu} b_{\lambda, \mu}^{\,\nu} L_\nu
    &
    \Delta(L_{\nu}) \ &=\ \sum_{\lambda,\mu} c_{\lambda, \mu}^{\,\nu} L_\lambda \otimes L_\mu
        \ =\ \sum_{\mu} L_{\nu / \mu} \otimes L_\mu \\
\label{eq: R structure}
    R_\lambda \cdot R_\mu \ &=\ \sum_{\nu} c_{\lambda, \mu}^{\,\nu} R_\nu
    &
    \Delta(R_{\nu}) \ &=\ \sum_{\lambda,\mu} b_{\lambda, \mu}^{\,\nu} R_\lambda \otimes R_\mu
        \ =\ \sum_{\mu} R_{\nu / \mu} \otimes R_\mu \,. \rule[0pt]{0pt}{16pt}
\end{align}
The \demph{skew elements} $L_{\nu / \mu}$ and $R_{\nu / \mu}$ defined above co-multiply according to
\begin{gather}
\label{eq: LR skew coproduct} \Delta(L_{\tau/\sigma}) \ =\
\sum_{\pi,\rho} c_{\pi,\rho,\sigma}^{\,\tau}\, L_\pi \otimes L_\rho
\qquad\qquad \Delta(R_{\tau/\sigma}) \ =\ \sum_{\pi,\rho}
b_{\pi,\rho,\sigma}^{\,\tau}\, R_\pi \otimes R_\rho \,.
\end{gather}
(Note that the structure of $H^*$ can be recovered from the structure of $H$. Thus, we may suppress the analogs of \eqref{eq: R structure} and the second formula in \eqref{eq: LR skew coproduct} in the coming sections.)

Finally, suppose that the antipode acts on $H$ in the $L$-basis according to the formula
\begin{gather}\label{eq: antipode}
S(L_\rho) \ =\ (-1)^{\mathsf e(\rho)} L_{\rho^{\mathsf T}}
\end{gather}
for some functions $\mathsf e\colon \mathcal P \to \mathbb N$ and
$(\cdot)^{\mathsf T}\colon \mathcal P \to \mathcal P$.
Then Lemma \ref{Lem:Hopf} takes the following form.

\begin{theorem}[Algebraic Littlewood-Richardson formula]\label{th: algebraic rule}
For any $\lambda, \mu,\sigma,\tau \in \mathcal P$, we have
 \begin{equation}\label{eq: alg rule}
   L_{\mu/\lambda}\cdot L_{\tau/\sigma}\ =\
   \sum_{\pi,\rho,\lambda^-\!,\mu^+} (-1)^{\Red{\mathsf e(\rho)}} \;
   c^{\,\tau}_{\pi,\rho,\sigma}\; b^{\,\lambda}_{\lambda^-\!,\Red{\rho^{\mathsf T}}} \;
   b^{\,\mu^+}_{\mu,\pi} \; \; L_{\mu^+/\lambda^-} \,.
 \end{equation}
Swapping $L\leftrightarrow R$ and $b\leftrightarrow c$ in \eqref{eq: alg rule} yields the analog for the skew elements $R_{\mu / \lambda}$ in $H^*$.
\end{theorem}

\begin{proof}
 The actions in~\eqref{Eq:harpoon} together with the second formulas for
 the coproducts in~\eqref{eq: L structure} and \eqref{eq: R structure} show that $R_\lambda \harpoon L_\mu= L_{\mu/\lambda}$ and $L_\lambda \harpoon R_\mu= R_{\mu/\lambda}$.
Now use~\eqref{Two} and~\eqref{eq: L structure}--\eqref{eq: antipode} to obtain
 \begin{align*}
 L_{\mu/\lambda}\cdot L_{\tau/\sigma}\ =\  (R_\lambda\harpoon L_\mu)\cdot L_{\tau/\sigma} \ &=\
   \sum_{\pi, \rho} (-1)^{\mathsf e(\rho)} \;
         c^{\,\tau}_{\pi,\rho,\sigma}  \;
         \bigl((L_{\rho^{\mathsf T}}\harpoon R_\lambda)\harpoon(L_\mu\cdot L_\pi)\bigr)\\
   &=\ \sum_{\pi,\rho,\mu^+} (-1)^{\mathsf e(\rho)} \;
         c^{\,\tau}_{\pi,\rho,\sigma}  \;
         b^{\,\mu^+}_{\mu,\pi} \;
         \bigl( R_{\lambda/\rho^{\mathsf T}}\harpoon L_{\mu^+}\bigr)\\
   &=\ \sum_{\pi,\rho,\lambda^-\!,\mu^+} (-1)^{\mathsf e(\rho)} \;
          c^{\,\tau}_{\pi,\rho,\sigma}  \;
          b^{\,\lambda}_{\lambda^-\!,\rho^{\mathsf T}}  \;
          b^{\,\mu^+}_{\mu,\pi} \;
          (R_{\lambda^-}\harpoon L_{\mu^+})\,.
 \end{align*}
This equals the right hand side of~\eqref{eq: alg rule}, since
$R_{\lambda^-}\harpoon L_{\mu^+} = L_{\mu^+ / \lambda^-}$.
\end{proof}

\begin{remark}
{\rm
The condition \eqref{eq: antipode} is highly restrictive. It implies that the antipode $S$, as a linear map, is conjugate to a signed permutation matrix. Nevertheless, it holds for the Hopf algebras we consider. More generally, it holds if either $H$ or $H^*$ is commutative, for then $S$ is an involution \cite[Cor. 1.5.12]{Mont}.
}
\end{remark}

\section{Skew Littlewood-Richardson rule for Schur functions}\label{sec: Schur s}

The commutative Hopf algebra $\Blue{\Lambda}$ of symmetric functions is graded
and self-dual under the Hall inner product 
$\Blue{\langle \cdot,\cdot\rangle}\colon \Lambda\otimes\Lambda\to\mathbb{Q}$. 
A systematic study of $\Lambda$ from a Hopf algebra perspective appears in \cite{Zel:1981}. 
We follow the definitions and notation in Chapter I of~\cite{Mac}.
The Schur basis of $\Lambda$ (indexed by partitions) is self-dual, so \eqref{eq: L structure} and \eqref{eq: LR skew coproduct} become
 \begin{gather}\label{Eq:LR_rule}
   s_\lambda \cdot s_\mu\ =\ \sum_\nu c^{\,\nu}_{\lambda,\mu} s_\nu
   \qquad\qquad
   \Delta(s_\nu)\ =\ \sum_{\lambda,\mu} c^{\,\nu}_{\lambda,\mu} s_\lambda\otimes s_\mu
   \ =\ \sum_{\mu}  s_{\nu/\mu} \otimes s_{\mu}
\\
\label{Eq:skew_coproduct}
   \Delta(s_{\tau/\sigma})\ =\
   \sum_{\pi,\rho} c^{\,\tau}_{\pi,\rho,\sigma}\; s_\pi\otimes s_\rho \,,
 \end{gather}
where the $c^{\,\nu}_{\lambda,\mu}$ are the \demph{Littlewood-Richardson coefficients} and the $s_{\nu/\mu}$ are the \demph{skew Schur functions} \cite[I.5]{Mac}.
Combinatorial expressions for the $c_{\lambda,\mu}^{\,\nu}$ and inner products $\langle s_{\mu/\lambda},s_{\tau/\sigma}\rangle$ are derived using the Hopf algebraic structure of $\Lambda$ in \cite{Zel:1981}. 
The coefficients \demph{$c^{\,\tau}_{\pi,\rho,\sigma}$} occur in the triple product $s_\pi\cdot s_\rho\cdot s_\sigma$,
\[
   c^{\,\tau}_{\pi,\rho,\sigma}
    \ =\ \langle s_\pi\cdot s_\rho\cdot s_\sigma, \, s_\tau \rangle
    \ =\ \langle s_\pi\cdot s_\rho,\, s_{\tau/\sigma} \rangle
    \ =\ \langle s_\pi\otimes s_\rho,\, \Delta(s_{\tau/\sigma}) \rangle \,.
\]
Write $\rho'$ for the conjugate (matrix-transpose) of $\rho$. Then the action of the antipode is
 \begin{equation}\label{eq: s-antipode}
   S(s_\rho)\ =\ (-1)^{|\rho|}s_{\rho'}\,,
 \end{equation}
 which is just a twisted form of the fundamental involution $\omega$
 that sends $s_\rho$ to $s_{\rho'}$.
 Indeed, the formula $\sum_{i+j=n}(-1)^ie_ih_j=\delta_{0,n}$ shows that~\eqref{eq: s-antipode}
 holds on the generators $\{h_n\mid n\geq 1\}$ of $\Lambda$.
 The validity of~\eqref{eq: s-antipode} follows as both $S$ and $\omega$ are algebra maps.

Since $c_{\lambda^-\!,\rho'}^{\,\lambda}=0$ unless $|\rho|=|\lambda / \lambda^-|$, we may write \eqref{eq: alg rule} as
 \begin{equation}\label{eq: Schur-s alg rule}
   s_{\mu/\lambda}\cdot s_{\tau/\sigma}\ =\
   \sum_{\pi,\rho,\lambda^-\!,\mu^+} (-1)^{|\lambda/\lambda^-|} \;
   c^{\,\tau}_{\pi,\rho,\sigma}\; c^{\,\lambda}_{\lambda^-\!,{\rho'}} \;
   c^{\,\mu^+}_{\mu,\pi} \; \; s_{\mu^+/\lambda^-} \;.
 \end{equation}

We next formulate a combinatorial version of~\eqref{eq: Schur-s alg rule}.
Given partitions $\rho$ and $\sigma$, form the skew shape \Blue{$\rho*\sigma$} by placing $\rho$
southwest of $\sigma$.
Thus,
\[
  \mbox{if\ }\
    \rho\ =\ \raisebox{-.45\height}{\includegraphics{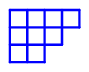}}
   \ \mbox{\ and\ }\
     \sigma\ =\ \raisebox{-.45\height}{\includegraphics{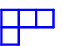}}\,
   \ \mbox{\ then\ }\
   \rho*\sigma\ =\ \raisebox{-.45\height}{\includegraphics{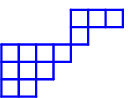}}\,.
\]
Similarly, if $R$ is a tableau of shape $\rho$ and $S$ a tableau of shape $\sigma$, then
$R*S$ is the skew tableau of shape $\rho*\sigma$ obtained by placing $R$ southwest of $S$.
Fix a tableau $T$ of shape $\tau$.
The Littlewood-Richardson coefficient $c^{\,\tau}_{\rho,\sigma}$ is the number of pairs $(R,S)$ of
tableaux of respective shapes $\rho$ and $\sigma$ with $R*S$ Knuth-equivalent to $T$.
See~\cite[Ch. 5, Cor. 2(v)]{Fu}.
Similarly, $c^{\,\tau}_{\pi,\rho,\sigma}$
is the number of triples $(P,R,S)$ of tableaux of respective shapes
$\pi$, $\rho$, and $\sigma$ with $P*R*S$ Knuth-equivalent to $T$.

Write \Blue{$\Sh(S)$} for the shape of a tableau $S$ and $S\equiv_K T$ if $S$ is Knuth-equivalent to $T$.

\begin{lemma}\label{L:Delta(skew)}
 Let $\sigma,\tau$ be partitions and fix a tableau $T$ of shape $\tau$.
 Then
\[
   \Delta (s_{\tau/\sigma})\ =\ \sum s_{\sh(R^-)} \otimes s_{\sh(R^+)}\,,
\]
the sum taken over triples $(R^-\!,R^+\!,S)$ of tableaux with $\Sh(S)=\sigma$ and
 $R^-*R^+*S \equiv_K T$. \qed
\end{lemma}

Note that $(\mu/\lambda)'=\mu'/\lambda'$ and the operation $*$ makes sense for
skew tableaux. If $S$ is a tableau of skew shape $\mu/\lambda$, put $|S|=|\mu/\lambda|=|\mu|-|\lambda|$.

\begin{theorem}[Skew Littlewood-Richardson rule]\label{th:two}
 Let $\lambda,\mu,\sigma,\tau$ be partitions and fix a tableau $T$ of shape $\tau$.
 Then
 \begin{equation}\label{eq:SLRII}
   s_{\mu/\lambda}\cdot s_{\tau/\sigma}\ =\
   \sum (-1)^{|S^-|} \; s_{\mu^+/\lambda^-}\ ,
 \end{equation}
the sum taken over triples $(S^-\!,S^+\!,S)$ of skew tableaux of respective shapes
 $(\lambda/\lambda^-)'$, $\mu^+/\mu$, and $\sigma$ such that
 $S^-*S^+*S\equiv_K T$.
\end{theorem}

\begin{remark}
{\rm
 If $T$ is the unique \demph{Yamanouchi tableau} of shape $\tau$ whose $i$th row contains only the letter $i$, then this is \Blue{{\it almost}} Conjecture 6.1 in~\cite{AM}.
Indeed, in this case $S$ is Yamanouchi of shape $\sigma$, so the sum is really over pairs of tableaux, and this explains the $\sigma$-Yamanouchi condition in \cite{AM}. The difference lies in the tableau $S^-$ and the reading word condition in~\cite{AM}.
It is an exercise in tableaux combinatorics that there is a bijection between the indices $(S^-\!,S^+)$ of Theorem~\ref{th:two} and the corresponding indices of Conjecture 6.1 in~\cite{AM}.
}
\end{remark}

\begin{proof}[Proof of Theorem \ref{th:two}]
 We reinterpret~\eqref{eq: Schur-s alg rule} in terms of tableaux.
Let $(R^-\!,R^+\!,S)$ be a triple of tableaux of partition shape with $\Sh(S)=\sigma$ and $R^-*R^+*S\equiv_K T$.
If $\Sh(R^-)=\rho$, then by~\cite[Ch. 5, Cor. 2(i)]{Fu},
 $c^{\,\lambda}_{\lambda^-\!,\rho'}=c^{\,\lambda'}_{(\lambda^-)',\rho}$ counts skew tableaux $S^-$ of shape $(\lambda/\lambda^-)'$ with $S^- \equiv_K R^-$.
 Likewise, if $\Sh(R^+)=\pi$, then $c^{\,\mu^+}_{\mu,\pi}$ counts skew tableaux $S^+$ of shape $\mu^+/\mu$ with $S^+\equiv_K R^+$.
Now~\eqref{eq: Schur-s alg rule} may be written as
\[
  s_{\mu/\lambda}\cdot s_{\tau/\sigma}\ =\
  \sum (-1)^{|S^-|} s_{\mu^+/\lambda^-}\,,
\]
summing over skew tableaux $(R^-\!,R^+\!,S^-\!,S^+\!,S)$ with
 $R^\pm$ of partition shape,
 $\Sh(S)=\sigma$,
 $R^-*R^+*S \equiv_K T$,
 $\Sh(S^+)=\mu^+/\mu$,
 $\Sh(S^-)=(\lambda/\lambda^-)'$,
and
 $S^\pm\equiv_K R^\pm$.

 Finally, note that $R^\pm$ is the unique tableau of partition shape Knuth-equivalent to $S^\pm$.
Since $S^-*S^+*S$ is Knuth-equivalent to $T$ (by transitivity of $\equiv_K$), we omit the unnecessary tableaux $R^\pm$ from the indices of summation and reach the statement of the theorem.
\end{proof}

%
\section{Skew Littlewood-Richardson rule for Schur $P$- and $Q$-functions}\label{sec: Schur PQ}

The self-dual Hopf algebra of symmetric functions has a natural self-dual subalgebra \Blue{$\Omega$}.
This has dual bases the Schur $P$- and $Q$-functions~\cite[III.8]{Mac}, which are indexed by \demph{strict partitions} $\lambda\colon \lambda_1>\dotsb>\lambda_l>0$.
Write \Blue{$\ell(\lambda)=l$} for the length of the partition $\lambda$.
As in Section \ref{sec: Schur s}, the constants and skew functions in the structure equations
\begin{gather}
\label{eq: Q structure}
    Q_\lambda\cdot Q_\mu \ =\ \sum_\nu g^{\,\nu}_{\lambda,\mu}\; Q_\nu
    \qquad\qquad
    \Delta(Q_\nu) \ =\ \sum_{\lambda,\mu} f^{\,\nu}_{\lambda,\mu} Q_\lambda\otimes Q_\mu
        \ =\ \sum_{\mu}  Q_{\nu/\mu}\otimes Q_{\mu}\, \\
\label{eq: PQ skew coproduct}
    \Delta(Q_{\tau/\sigma}) \ =\ \sum_{\pi,\rho} f_{\pi,\rho,\sigma}^{\,\tau}\, Q_\pi \otimes Q_\rho
\end{gather}
have combinatorial interpretations (see below).
Also, each basis $\{P_\lambda\}$ and $\{Q_\lambda\}$ is \emph{almost} self-dual in that $P_\lambda= 2^{-\ell(\lambda)}Q_\lambda$ and $g^{\,\nu}_{\lambda,\mu}=2^{\ell(\lambda)+\ell(\mu)-\ell(\nu)}f_{\lambda,\mu}^{\,\nu}$.

The algebra $\Omega$ is generated by the special $Q$-functions
$q_n=Q_{(n)}:=\sum_{i+j=n} h_ie_j$ \cite[III, (8.1)]{Mac}.
This implies that $S(q_n)=(-1)^n q_n$, from which we deduce that
\begin{align*}
S(Q_\rho)=(-1)^{|\rho|}Q_\rho \,.
\end{align*}
As $f_{\lambda^-\!,\rho}^{\,\lambda}=0$ unless $|\rho|=|\lambda /
\lambda^-|$, we may write the algebraic rule~\eqref{eq: alg rule} as
\begin{align}\label{eq: Schur-Q alg rule}
   Q_{\mu/\lambda}\cdot Q_{\tau/\sigma}\ =\
   \sum_{\pi,\rho,\lambda^-\!,\mu^+} (-1)^{|\lambda/\lambda^-|}
  \; f^{\,\tau}_{\pi,\rho,\sigma}\;
     g^{\,\lambda}_{\lambda^-\!,\rho} \; g^{\,\mu^+}_{\mu,\pi} \; \; Q_{\mu^+/\lambda^-}\;,
\end{align}
with a similar identity holding for $P_{\mu/\lambda}\cdot P_{\tau/\sigma}$ (swapping $P\leftrightarrow Q$ and $f\leftrightarrow g$).

We formulate two combinatorial versions of~\eqref{eq: Schur-Q alg rule}.
Strict partitions $\lambda,\mu$ are written as \demph{shifted Young diagrams} (where row $i$ begins in column $i$).
Skew shifted shapes $\lambda/\mu$ are defined in the obvious manner:
\[
  \mbox{if\ }\
    \lambda\ =\ 431 \ =\ \raisebox{-.45\height}{\includegraphics{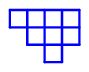}}
   \ \mbox{\ and\ }\
     \mu\ =\ 31 \ =\ \raisebox{-.4\height}{\includegraphics{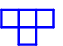}}\,,
   \ \mbox{\ then\ }\
   \lambda/\mu\ =\ \raisebox{-.45\height}{\includegraphics{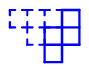}}\ =\ \raisebox{-.45\height}{\includegraphics{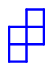}}\,.
\]
In what follows, tableaux means \emph{semi-standard (skew) shifted tableaux on a marked alphabet}~\cite[III.8]{Mac}.
We use shifted versions of the jeu-de-taquin and plactic equivalence from~\cite{Sag} and~\cite{Ser}, denoting the corresponding relations by \demph{$\eqsj$} and \demph{$\eqsp$}, respectively.
Given tableaux $R,S,T,$ we write $R*S\eqsp T$ when representative words $u,v,w$ (built via ``$\mathsf{mread}$''~\cite[\S2]{Ser}) of the corresponding shifted plactic classes satisfy $uv\eqsp w$.

Stembridge notes (following \cite[Prop. 8.2]{Stem}) that for a fixed tableau $M$ of shape $\mu$,
\begin{gather}\label{eq: Stem LR}
    f_{\lambda,\mu}^{\,\nu} \ =\ \#\bigl\{ \hbox{skew tableaux }L :
    \Sh(L)=\nu/\lambda \hbox{ and } L\eqsj M \bigr\}.
\end{gather}
Serrano has a similar description of these coefficients in terms of $\eqsp$. Fixing a tableau $T$ of shape $\tau$, it follows from \cite[Cor. 1.15]{Ser} that
the coefficient $f_{\pi,\rho,\sigma}^{\,\tau}$ in $P_\pi \cdot P_\rho \cdot P_\sigma = \sum_{\tau} f_{\pi,\rho,\sigma}^{\,\tau} P_\tau$ is given by
\begin{gather}\label{eq: Ser triple LR}
    f_{\pi,\rho,\sigma}^{\,\tau} \ =\ \#\bigl\{ (P,R,S) :
    \Sh(P)=\pi,\, \Sh(R)=\rho,\, \Sh(S)=\sigma,\hbox{ and }
    P*R*S \eqsp T \bigr\} .
\end{gather}

If $T$ is a tableau of shape $\lambda$, write $\ell(T)$ for $\ell(\lambda)$.
The formula relating the $g$'s and $f$'s combines with \eqref{eq: Stem LR} and \eqref{eq: Ser triple LR} to give our next result.

\begin{theorem}[Skew Littlewood-Richardson rule]\label{th: Schur-Q}
 Let $\lambda,\mu,\sigma,\tau$ be shifted partitions and fix a tableau $T$ of shape $\tau$.
 Then
 \begin{equation}\label{eq: Schur-Q SLR}
   Q_{\mu/\lambda}\cdot Q_{\tau/\sigma}\ =\
   \sum (-1)^{|\lambda/\lambda^-|}\, 2^{\ell(R^-)+\ell(R^+)+\ell(\lambda^-)+\ell(\mu)-\ell(\lambda)-\ell(\mu^+)} \;
   Q_{\mu^+/\lambda^-}\ ,
 \end{equation}
the sum taken over quintuples $(R^-\!,R^+\!,S^-\!,S^+\!,S)$ with
 $R^\pm$ of partition shape,
 $\Sh(S)=\sigma$,
 $R^-*R^+*S \eqsp T$,
 $\Sh(S^+)=\mu^+/\mu$,
 $\Sh(S^-)=(\lambda/\lambda^-)$,
and
 $S^\pm\eqsj R^\pm$. \qed
\end{theorem}

Serrano conjectures an elegant combinatorial description \cite[Conj.
2.12 and Cor. 2.13]{Ser} of the structure constants
$g_{\lambda,\mu}^{\,\nu}$ in \eqref{eq: Q structure}: For any
tableau $M$ of shape $\mu$, he conjectured
\begin{gather}\label{eq: Ser LR}
    g_{\lambda,\mu}^{\,\nu} \ =\ \#\bigl\{ \hbox{skew tableaux }L :
    \Sh(L) = \nu/\lambda \hbox{ and } L \eqsp M \bigr\}.
\end{gather}
(Note that if $S,T$ are tableaux, then $S\eqsp T$ does not necessarily imply that $S\eqsj T$.)
This leads to a conjectural reformulation of Theorem \ref{th: Schur-Q} in the spirit of Theorem \ref{th:two}.

\begin{conjecture}[Conjectural Skew Littlewood-Richardson rule]\label{th: conj Schur-Q}
 Let $\lambda,\mu,\sigma,\tau$ be partitions and fix a tableau $T$ of shape $\tau$.
 Then
 \begin{equation}\label{eq: conj Schur-Q SLR}
   Q_{\mu/\lambda}\cdot Q_{\tau/\sigma}\ =\
   \sum (-1)^{|S^-|} \;
   Q_{\mu^+/\lambda^-}\ ,
 \end{equation}
the sum taken over triples $(S^-\!,S^+\!,S)$ of skew tableaux of respective shapes
$(\lambda/\lambda^-)$, $\mu^+/\mu$, and $\sigma$ such that
${S^-}*{S^+}*{S} \eqsp T$.
\end{conjecture}

\begin{proof}
There is a unique shifted tableau $R$ in any shifted plactic class~\cite[Thm. 2.8]{Ser}.
So the conditions $S^\pm \eqsp R^\pm$ and ${R^-}*{R^+}*{S} \eqsp T$ in \eqref{eq: Ser triple LR} and \eqref{eq: Ser LR} may be replaced with the single condition  ${S^-}*{S^+}*{S} \eqsp T$.
\end{proof}

%
\section{Skew Littlewood-Richardson rule for noncommutative Schur functions}\label{sec: Schur NSym}

The Hopf algebra of noncommutative symmetric functions was
introduced, independently, in \cite{GKLLRT:1995,MalReu:1995} as the (graded) dual to the commutative Hopf algebra of quasisymmetric functions. We consider the dual bases (indexed by compositions)
$\{F_\alpha\}$ of Gessel's quasisymmetric functions and $\{R_\alpha\}$ of
noncommutative ribbon Schur functions.
The structure constants in
\begin{gather}
 \label{eq: FR structure}
     R_\alpha\cdot R_\beta \ =\ \sum_\gamma b^{\,\gamma}_{\alpha,\beta}\; R_\gamma
     \qquad\qquad
     \Delta(R_\gamma) \ =\ \sum_{\alpha,\beta} c^{\,\gamma}_{\alpha,\beta} R_\alpha\otimes R_\beta
        \ =\ \sum_{\beta}  R_{\gamma/\beta}\otimes R_{\beta}\, \\
\label{eq: FR skew coproduct}
    \Delta(R_{\tau/\sigma}) \ =\ \sum_{\pi,\rho,\sigma} c_{\pi,\rho,\sigma}^{\,\tau}\, R_\pi \otimes R_\rho \,
\end{gather}
may be given combinatorial meaning via the descent map $\des\colon \mathfrak{S}_n \to \Gamma_n$ from permutations to compositions and a section of it $\word\colon \Gamma_n \to \mathfrak{S}_n$.
These maps are linked via \demph{ribbon diagrams}, edge-connected skew Young diagrams (written in the french style), with no $2\times 2$ subdiagram present.
By way of example,
\begin{align*}
  \des\colon\ {148623795} &\rightarrow
    \raisebox{-.35\height}{\includegraphics{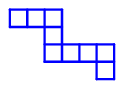}} \rightarrow
    {3141}, &
  \word\colon\ {3141} &\rightarrow
    \raisebox{-.35\height}{\includegraphics{ribbon.eps}} \rightarrow
    {789623451}.
\end{align*}
(In the intermediate step for $\des(w)$, new rows in the ribbon begin at {descents} of $w$.
In the intermediate step for $\word(\alpha)$, the boxes in the ribbon are filled left-to-right, bottom-to-top.)

A ribbon $\alpha$ may be extended by a ribbon $\beta$ in two ways: affixing $\beta$ to the rightmost edge or bottommost edge of $\alpha$ (written $\alpha\ljoin \beta$ and $\alpha\ujoin \beta$, respectively):
\begin{align*}
  {311 \ljoin{\scriptscriptstyle\,}31}\colon &
    \raisebox{-.35\height}{\includegraphics{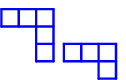}} \rightarrow
    \raisebox{-.35\height}{\includegraphics{ribbon.eps}}\,, &
  {31\ujoin{\scriptscriptstyle\,}41}\colon &
    \raisebox{-.35\height}{\includegraphics{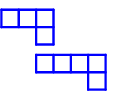}} \rightarrow
    \raisebox{-.35\height}{\includegraphics{ribbon.eps}}\,.
\end{align*}
If a ribbon $\gamma$ is formed from $\alpha$ and $\beta$ in either of these two ways, we write $\gamma\in \alpha\join \beta$.
The coefficient $b_{\alpha,\beta}^{\,\gamma}$ is $1$, if $\gamma\in \alpha\join \beta$, and $0$ otherwise.
If $\ast$ is the shifted shuffle product on permutations (see (3.4) in \cite{MalReu:1995}), then the coefficient $c_{\alpha,\beta}^{\,\gamma}$ is the number of words $w$ in $\word(\alpha)\ast\word(\beta)$ such that $\des(w) = \gamma$. The coefficient $c_{\pi,\rho,\sigma}^{\,\tau}$ has the analogous description.

Antipode formulas for the distinquished bases were found, independently, in \cite{Ehr:1996,MalReu:1995}:
\begin{equation*}
S(F_\alpha) \ =\ (-1)^{|\alpha|}F_{\tild \alpha} \qquad\hbox{and}\qquad
S(R_\alpha) \ =\ (-1)^{|\alpha|}R_{\tild \alpha} \,,
\end{equation*}
where $\tild \alpha$ is the conjugate of $\alpha$ (in the sense of french style skew partitions).
For example, $\tild{(3141)} = 211311$. The descriptions of the antipode and structure constants in \eqref{eq: FR structure} and \eqref{eq: FR skew coproduct} give a formula for the product of two skew ribbon Schur functions.

\begin{theorem}[Skew Littlewood-Richardson rule]\label{th: Schur-R}
Let $\alpha,\beta,\sigma,\tau$ be compositions. 
Then
 \begin{equation*}\label{eq: Schur-R SLR}
   R_{\beta/\alpha}\cdot R_{\tau/\sigma}\ =\
   \sum (-1)^{|{\rho}|} \;
   R_{\beta^+/\alpha^-}\ ,
 \end{equation*}
the sum taken over factorizations $\alpha\in \alpha^-\join \tild{\rho}$,
extensions $\beta^+\in \beta\join \pi$, and words $w$ in the
shuffle product $\word(\pi)\ast \word({\rho}) \ast \word(\sigma)$ such that $\des(w) = \tau$. \hfill \qed
\end{theorem}

\begin{remark}\label{rem: F skew}
{\rm
The nonzero skew ribbon Schur functions $R_{\beta/\alpha}$ do not correspond to skew ribbon shapes in a simple way.
For example, $111$ is not a (connected) sub-ribbon of $221$, yet $R_{221/111} = R_{2} + R_{11}\neq0$.
Contrast this with the skew functions $F_{\beta/\alpha}$, where
$F_{\beta/\alpha}\neq0$ if and only if $\beta\in \omega\join \alpha$ for some ribbon $\omega$. That is, $F_{\beta/\alpha}=F_\omega$. Thus we may view $\alpha$ as the last $|\alpha|$ boxes of the ribbon $\beta$ and $\beta/\alpha$ as the complementary ribbon $\omega$. Interpreting $F_{\beta/\alpha}\cdot F_{\tau/\sigma}$ alternately as a product of ordinary functions or skew functions yields the curious identity
\begin{equation}\label{eq: Gessel-F alg rule}
F_{\beta/\alpha}\cdot F_{\tau/\sigma} \ =\
   \sum_{\gamma} c_{\beta/\alpha,\tau/\sigma}^{\,\,\gamma} \; F_\gamma \ =\
   \sum_{\pi,\rho,\alpha^-\!,\beta^+} (-1)^{|\rho|}
  \; b^{\,\tau/\sigma}_{\pi,\rho}\;
     c^{\,\alpha}_{\alpha^-\!,\tild{\rho}} \; c^{\,\beta^+}_{\beta,\pi} \; \; F_{\beta^+/\alpha^-}\; .
\end{equation}
}
\end{remark}

%
\section{Skew $k$-Pieri rule for $k$-Schur functions}\label{sec:kschur}
Fix an integer $k \geq 1$. Let $\Lambda_{(k)}$ denote the Hopf
subalgebra of the Hopf algebra of symmetric functions generated by
the homogeneous symmetric functions $h_1,h_2,\ldots,h_k$.  Let
$\Lambda^{(k)}$ denote the Hopf-dual quotient Hopf algebra of
symmetric functions.  We consider the dual bases
$\{s_\lambda^{(k)}\} \subset \Lambda_{(k)}$ and $\{F_\lambda^{(k)}\}
\subset \Lambda^{(k)}$ of $k$-Schur functions and dual $k$-Schur
functions of \cite{LLMS,LM}, also called strong Schur functions and
weak Schur functions in \cite{LLMS}. The $k$-Schur functions were
first introduced by Lapointe, Lascoux, and Morse in the context of
Macdonald polynomials, and were later shown by Lam to represent
Schubert classes in the affine Grassmannian of $\mathrm{SL}(k{+}1,\mathbb C)$. We refer
the reader to the references in \cite{LLMS}.

Here $\lambda$ varies over all $k$-bounded partitions, that is,
those partitions satisfying $\lambda_1 \leq k$. There is an
involution $\lambda \mapsto \lambda^{\omega_k}$ on $k$-bounded
partitions called \demph{$k$-conjugation}.  We have
\begin{equation*}
  S(s_\lambda^{(k)}) \ =\ (-1)^{|\lambda|} s_{\lambda^{\omega_k}}^{(k)}
\qquad \hbox{and} \qquad
  S(F_\lambda^{(k)}) \ =\ (-1)^{|\lambda|} F_{\lambda^{\omega_k}}^{(k)}.
\end{equation*}

If $\lambda = (r)$ is a one-part partition, then $s_\lambda^{(k)}
= h_r$ is a homogeneous symmetric function.  We have the
$k$-Pieri and dual $k$-Pieri rules \cite{LLMS,LM} (called weak
and strong Pieri rules in \cite{LLMS})
\begin{equation}
\label{eq:kPieri}
  s_\lambda^{(k)} \cdot h_r \ =\
    \sum_{\lambda \warrow{r} \mu } s_\mu^{(k)}
\qquad \hbox{and} \qquad
  F_\lambda^{(k)} \cdot h_r \ =\
    \sum_{\lambda \sarrow{r} \mu} F_\mu^{(k)}
\end{equation}
for $r \leq k$.  Here $\lambda \warrow{r} \mu$ denotes an \demph{$r$-weak strip} connecting $\lambda$ and $\mu$---present if and only if both
$\mu/\lambda$ and $\mu^{\omega_k}/\lambda^{\omega_k}$ are horizontal
$r$-strips. The notation $\lambda \sarrow{r} \mu$ denotes an \demph{$r$-strong
strip} as introduced in \cite{LLMS}, which we will not define here.
The terminology comes from the relationship with the weak and strong
(Bruhat) orders of the affine symmetric group. We remark that there
may be distinct strong strips $\lambda \sarrow{r} \mu$ and $(\lambda
\sarrow{r} \mu)'$ which start and end at the same partition, so that
the second Pieri rule of \eqref{eq:kPieri} may have multiplicities.
(Strictly speaking, the strong strips in \cite{LLMS} are built on
$(k+1)$-cores, and our $\lambda \sarrow{r} \mu$ denotes the strips
obtained after applying a bijection between $(k+1)$-cores and
$k$-bounded partitions.)

We define skew functions $s_{\lambda/\mu}^{(k)}$ and $F_{\lambda/\mu}^{(k)}$ using
\eqref{eq: L structure} and \eqref{eq: R structure}.  There is an
explicit combinatorial description of $F_{\lambda/\mu}^{(k)}$ in
terms of the weak tableaux of \cite{LLMS}, but only a conjectured
combinatorial description of $s_{\lambda/\mu}^{(k)}$
\cite[Conj. 4.18(3)]{LLMS}.

\begin{theorem}[Skew $k$-Pieri (or weak Pieri) rule]
\label{thm:skewkPieri} For $k$-bounded partitions $\lambda, \mu$,
and $r \leq k$,
\begin{equation*}
  s_{\mu/\lambda}^{(k)} \cdot h_r \  =\ \sum_{i+j=r} (-1)^j \,
  \sum_{\substack{\mu \warrow{i} \mu^+ \\[.25ex]
  (\lambda^-)^{\omega_k} \warrow{j} \lambda^{\omega_k} }}
  s_{\mu^+/\lambda^-}^{(k)}
\end{equation*}
\end{theorem}
\begin{proof}
In Theorem \ref{th: algebraic rule}, take $L_{\tau/\sigma} = h_r$.  For
$c^{\,\tau}_{\pi,\rho,\sigma}$, use the formula $\Delta(h_r) =
\sum_{i+j=r} h_{i} \otimes h_{j}$, and for
$b^{\,\lambda}_{\lambda^-\!,\rho^{\,\omega_k}}$ and $b^{\,\mu^+}_{\mu,\pi}$,
use \eqref{eq:kPieri}.
\end{proof}

\begin{theorem}[Skew dual $k$-Pieri (or strong Pieri) rule]
\label{thm:skewdualkPieri} For $k$-bounded partitions $\lambda,
\mu$, and $r \leq k$,
\begin{equation*}
  F_{\lambda/\mu}^{(k)} \cdot h_r \ =\ \sum_{i+j=r} (-1)^j \,
  \sum_{\substack{\lambda \sarrow{i} \lambda^+  \\[.25ex]
  (\mu^-)^{\omega_k} \sarrow{j} \mu^{\omega_k} }} F_{\lambda^+/\mu^-}^{(k)}
\end{equation*}
\begin{proof}
Identical to the proof of Theorem \ref{thm:skewkPieri}.
\end{proof}

\end{theorem}

As an example, let $k = 2$, $r = 2$, $\mu = (2,1,1)$, and $\lambda =
(1)$.  Then Theorem \ref{thm:skewkPieri} states that
\begin{equation*}
  s^{(2)}_{211/1} \cdot h_2 \ =\ s^{(2)}_{2211/1} - s^{(2)}_{2111} \,,
\end{equation*}
which one can verify using \eqref{eq:kPieri} and the expansions
$s^{(2)}_{211/1} = s^{(2)}_{21}+s^{(2)}_{111}$ and $s^{(2)}_{2211/1}
= 2\,s^{(2)}_{2111} + s^{(2)}_{221}$.  Theorem
\ref{thm:skewdualkPieri} states that
\begin{equation*}
  F^{(2)}_{211/1} \cdot h_2 \ =\ 3\, F^{(2)}_{222/1} + 5\,F^{(2)}_{2211/1} + 3\,F^{(2)}_{21111/1}
  + 3\, F^{(2)}_{111111/1} -2\,F^{(2)}_{221} - 3\,F^{(2)}_{2111} - 2\, F^{(2)}_{11111}.
\end{equation*}
One can verify that both sides are equal to
$6\,F^{(2)}_{221} + 5\,F^{(2)}_{2111} + 4\,F^{(2)}_{11111}$.

\section{Skew Pieri rule for affine Grassmannian of the symplectic group}

Fix $n\geq 1$.  The Hopf algebra $\Omega$ of Section \ref{sec: Schur
PQ} contains a Hopf subalgebra $\Omega_{(n)}$ generated by the Schur
$P$-functions $P_1,P_3,\ldots,P_{2n-1}$.  In \cite{LSS}, it was
shown that $\Omega_{(n)}$ is isomorphic to the homology ring
$H_*(\Gr_{\mathrm{Sp}(2n,\C)})$ of the affine Grassmannian of the
symplectic group $\mathrm{Sp}(2n,\C)$. A distinguished basis
$\{P_w^{(n)}\} \subset \Omega_{(n)}$, representing the Schubert
basis, was studied there.  The symmetric functions $P_w^{(n)}$ are
shifted versions of the $k$-Schur functions of Section
\ref{sec:kschur}.

The indexing set for the basis $\{P_w^{(n)}\}$ is the set
$\widetilde{C}^0_n$ of affine Grassmannian type $C$ permutations: they
are the minimal length coset representatives of $C_n$ in
$\widetilde{C}_n$. A lower Bruhat order ideal $\Z \subset \widetilde{C}_n$
of the affine type $C$ Weyl group is defined in \cite{LSS}.
%
Let $\Z_j\subset \Z$ denote those $v\in \Z$ with length $\ell(v)=j$.
%
For each $v \in \Z$, there is a nonnegative integer $\Blue{c(v)} \in
{\mathbb Z}_{\geq 0}$, called the \demph{number of components} of $v$.  We
note that $c({\rm id}) = 0$. With this notation, for each $1 \leq j
\leq 2n-1$, we have the Pieri rule \cite[Thms. 1.3 and 1.4]{LSS}
\begin{equation}
\label{eq:LSSPieri}
  P_w^{(n)} \cdot P_j \ =\ \sum_{v \in \Z_j} 2^{c(v) - 1}\; P_{vw}^{(n)}\, ,
\end{equation}
where the sum is over all $v \in \Z_j$ such that
$vw \in \widetilde{C}^0_n$, and
$\ell(vw) = \ell(v) + \ell(w)$.

It follows from the discussion in Section \ref{sec: Schur PQ} that
the antipode acts on the $P_j$ by $S(P_j) = (-1)^{j} P_j$.  We
define $P_{w/v}^{(n)}$ using \eqref{eq: L structure}.

\begin{theorem}[Skew Pieri rule]
For $w, v \in \widetilde{C}_n^0$, and $r \leq 2n-1$,
\begin{equation*}
  P_{w/v}^{(n)} \cdot P_r \ =\ \sum_{i+j=r} (-1)^j \,
  \sum_{\substack{u \in \Z_{i}\\ z \in \Z_{j}}} 2^{c(u) + c(z) - 1} \; P_{uw/z^{-1}v}^{(n)}\, ,
\end{equation*}
where the sum is over all $u \in \Z_{i}$ and $z \in \Z_{j}$ such that
$uw, z^{-1}v \in \widetilde{C}^0_n$,
$\ell(uw) = \ell(u) + \ell(w)$, and
$\ell(z^{-1}v) + \ell(z) = \ell(v)$.
\end{theorem}

\begin{proof}
In Theorem \ref{th: algebraic rule}, take $L_{\tau/\sigma}=P_r$ and use \eqref{eq:LSSPieri}. For the constants $c_{\pi,\rho,\sigma}^{\,\tau}$, use the formula
\begin{equation*}
  \Delta(P_r) \ =\
  1 \otimes P_r + P_r \otimes 1 + 2 \; \sum_{0 < j < r} P_{r-j} \otimes P_{j}\,.
\end{equation*}
If $0 < j < r$, the product
$c_{\pi,\rho,\sigma}^{\,\tau} b_{\lambda^-\!,\rho'}^{\,\lambda} b_{\mu,\pi}^{\,\mu^+}$ in \eqref{eq: alg rule} becomes $2 \cdot 2^{c(u)-1} \cdot 2^{c(z) - 1} = 2^{c(u)+c(z)-1}$.
If $j = 0$ (resp., $j=r$), it becomes $1\cdot 2^{c(u)-1}\cdot1 = 2^{c(u)+c(z)-1} $
(resp., $1\cdot1\cdot 2^{c(z)-1} =2^{c(u)+c(z)-1}$).
\end{proof}


\providecommand{\bysame}{\leavevmode\hbox to3em{\hrulefill}\thinspace}
\providecommand{\MR}{\relax\ifhmode\unskip\space\fi MR }
\providecommand{\MRhref}[2]{%
  \href{http://www.ams.org/mathscinet-getitem?mr=#1}{#2}
}
\providecommand{\href}[2]{#2}

\end{document}